\documentclass[10pt,twoside]{article}
\usepackage{amsfonts}
\usepackage{fancyhdr}
\usepackage{titlesec}
\usepackage{cite}
\usepackage{ifthen}
\usepackage{pifont}
\usepackage{stmaryrd}
\usepackage{setspace}
\usepackage{indentfirst}
\usepackage{graphicx}
\usepackage{amsmath,amssymb,amscd,bbm,amsthm,mathrsfs,dsfont}
\input amssym.def

\titleformat{\section}{\large\bfseries}{\arabic{section}}{1em}{}
\newboolean{first}
\setboolean{first}{true}

\newtheorem{theorem}{Theorem}[section]
\newtheorem{lemma}{Lemma}[section]
\newtheorem{rem}{Remark}[section]

\catcode`@=11
\@addtoreset{equation}{section}

\allowdisplaybreaks[2]

\begin{document}

\setlength\abovedisplayskip{2pt}
\setlength\abovedisplayshortskip{0pt}
\setlength\belowdisplayskip{2pt}
\setlength\belowdisplayshortskip{0pt}

\title{\bf \Large Entire Solutions of the Systems of Complex Differential-Difference Equations
 \author{  {\bf Xianfeng Su$^{1~2}$ ~~ Qingcai  Zhang$^{1~*}$ } \\
 {\small 1 School of Information, Renmin University  of  China,   \,Beijing, 100872,  China}\\
 {\small 2 School of Information, Huaibei Normal University,   \,Huaibei, 235000,  China}\\
{\small
E-mail: suxf2006@ruc.edu.cn ~~zhangqcrd@163.com
 }
 }
 \date{}} \maketitle
\footnote{$^{*}$Corresponding author. }
\footnote{This work is  supported by the National Natural Science Foundation of China (11171013), the Natural Science Foundation of the
 Education Department of Anhui Province of China (KJ2015A323) and the outstanding young talents program of the
 Education Department of Anhui Province of China(gxyq2017153).}

\begin{center}
\begin{minipage}{135mm}

{\bf \small Abstract}\hskip 2mm {\small  Using Nevanlinna's value
distribution theory and the complex difference theory , the existence of the finite order of  the entire solutions of the Fermat type of the complex differential-difference equation and the systems of complex differential-difference equations  are investigated, and obtained some results, which extends  Liu and Gao's results.}


{\bf \small 2010 MR Subject Classification } \hskip 2mm {\small 30D35; 39A10. }

\end{minipage}
\end{center}

\thispagestyle{fancyplain} \fancyhead{}
\fancyhead[L]{\textit{}\\
} \fancyfoot{} \vskip 10mm

\section{Introduction and notation}    
In this paper, we use the standard notations of the of Nevanlinna's value distribution theory, such as the characteristic function $T(r,w)$, proximity function $m(r,w)$, counting function $N(r,w)$ etc. and complex difference theory(see \cite{s1}-\cite{s2}). The order of the meromorphic function $w$ is given by
$$\rho(w)=\limsup_{r\rightarrow \infty}\frac{\log^{+} T(r,w)}{\log r}.$$
If $T(r,g)=S(r,w)$, we call $g(z)$ is a small function of $w(z)$, where $S(r,w)=o(T(r,w))$~$(r\rightarrow \infty )$, outside of a possible exceptional set with finite logarithmic measure.
\par Many authors, such as  Gross \cite{s3}, Yang\cite{s4}, Yang-Li \cite{s5}, Yi-Yang \cite{s6} etc. have investigated related problem of the Fermat type of equation, and obtained a series of important results. Recently, the investigation of complex difference has attracts a lot of attention and some significate results have been obtained, such as Chiang-Feng \cite{s7}, as well as Halburd-Korhonen \cite{s8}, have both independently established a difference analog of the logarithmic derivative lemma. Some authors considered the complex Fermat type of differential-difference equation and the systems of complex differential-difference equation, and also gained some interesting results(see\cite{s9}-\cite{s12}).
\par  The remainder part  of this paper is organized as follows: In Section 2, we investigate entire solutions  of the Fermat type of complex differential-difference equations, and gain some results, which generalize and improve the result of \cite{s9}. In section 3, we study the entire solutions of the Fermat type systems of  complex differential-difference equations. our results extend Gao's results \cite{s11}. In Section 4, we introduce a series  of  lemmas. They are essential to the proof of theorems. In Section 5, we prove theorems for Fermat type of complex differential-difference equations by lemmas given in Section 4. In Section 6, we mainly confirm our results for the Fermat type of systems of complex difference equations by our lemmas.
\section{Results on the Fermat type of complex differential-difference equations}
The following equation
\begin{align}\label{E:2.1} f(z)^{n}+g(z)^{n}=1
\end{align}
can is regard as the Fermat Diophantine equations $x^{n}+y^{n}=1$ ($n$ is a positive integer) over function fields. Gross \cite{s3} proved that Eq.\eqref{E:2.1} has no transcendental meromorphic solutions for $n \geq 4$, and has no entire solutions for $n \geq 3$. Yang \cite{s4} considered the following the Fermat type functional equations :
\begin{align}\label{E:2.2} f(z)^{n}+g(z)^{m}=1,
\end{align}
and got the following result.\\
{\bf Theorem A. (see \cite{s4})} Let $m$, $n$ be positive integers satisfying $\frac{1}{n}+\frac{1}{m}<1$. Then then equation \eqref{E:2.2} does not have  non-constant entire solutions $f(z)$ and $g(z)$.
\par In 2012, Liu-Cao-Cao\cite{s9} studied the Fermat type of complex differential-difference equations as below:
\begin{align}\label{E:2.31}
w^{\prime}(z)^{n}+w(z+c)^{m}=1,~~~~~~~~~~
\end{align}
\begin{align}\label{E:2.41}
w^{\prime}(z)^{n}+[w(z+c)-w(z)]^{m}=1,
\end{align}
where $m$, $n$ are positive integers.
They gained the following four theorems:\\
{\bf Theorem B.(\cite{s9})} Let $m$, $n$ be positive integers and $m\neq n$. Then Eq.\eqref{E:2.31} has no transcendental entire solutions with finite order.\\
{\bf Theorem D.(\cite{s9})} Suppose $m$, $n$ are positive integers satisfying $m\neq n$ and $\max\{m,n\}>1$. Then Eq.\eqref{E:2.41} has no transcendental entire solutions with finite order.
\par Liu-Cao-Cao\cite{s9} also considered the special cases $m=n=2$, got the following results.
\begin{align}\label{E:2.3}
w^{\prime}(z)^{2}+w(z+c)^{2}=1,~~~~~~~~~~~~
\end{align}
\begin{align}\label{E:2.4}
w^{\prime}(z)^{2}+[w(z+c)-w(z)]^{2}=1.
\end{align}
{\bf Theorem C.(\cite{s9})}  The transcendental entire solutions with finite order of Eq.\eqref{E:2.3} must satisfy $w(z)=\sin{(z+B)}$, where $B$ is a constant  and $c=k\pi$, $k$ can be any integer.\\
{\bf Theorem E.(\cite{s9})}  The transcendental entire solutions with finite order of Eq.\eqref{E:2.4} must satisfy $w(z)=\frac{1}{2}\sin{(2z+B)}$, where $B$ is a constant  and $c=k\pi+\frac{\pi}{2}$, $k$ can be  any integer.
\\
\par The first aim of this paper is to generalize the results in \cite{s9}. We mainly consider the following the Fermat type complex differential-difference equations:
\begin{align}\label{E:2.5}
w^{\prime}(z)^{n}+w(qz+c)^{m}=1,~~~~~~~~~~
\end{align}
\begin{align}\label{E:2.6}
w^{\prime}(z)^{n}+[w(qz+c)-w(z)]^{m}=1,
\end{align}
where $q(\neq0)$ and $c$ are complex constants.
Our results can be stated as the following five theorems.
\begin{theorem}  Let $m\geq2$, $n\geq2$ be positive integers and at lest one be not equal to 2. Then Eq.\eqref{E:2.5} and  Eq. \eqref{E:2.6} do not have non-constant entire solutions.
\end{theorem}
Next, we consider the case: $m=n=2$. That is
\begin{align}\label{E:2.7}
w^{\prime}(z)^{2}+w(qz+c)^{2}=1,~~~~~~~~~~~
\end{align}
\begin{align}\label{E:2.8}
w^{\prime}(z)^{2}+[w(qz+c)-w(z)]^{2}=1,
\end{align}
where $q(\neq0)$ and $c$ are complex constants. We have
\begin{theorem} According to the different value of $q$, we have
\par (1) If $q\neq \pm 1$, then  Eq.\eqref{E:2.7} does not posses any  finite order of transcendental entire solutions.
\par (2) If $q=1$, then Eq.\eqref{E:2.7} has the finite order of transcendental entire solution $$w(z)=\sin{(z+B)},$$ where $B$ is a constant and $c=k\pi$ $(k\in \mathds{Z})$.
\par (3) If $q=-1$, then Eq.\eqref{E:2.7} has the finite order of transcendental entire solution $$w(z)=\sin(z+B),$$ where $B$ is a constant  and $2B+c=k\pi$ $(k\in \mathds{Z})$.
\end{theorem}

Obviously, when $q=1$, then Eq.\eqref{E:2.7} reduce to Eq.\eqref{E:2.3}

\begin{theorem} According to the different value of $q$, we get
\par(1) If $q\neq \pm 1$, then  Eq.\eqref{E:2.8} does not posses any  finite order of transcendental entire solution.
\par (2) If $q=1$, then Eq.\eqref{E:2.8} has the finite order of transcendental entire solution $$w(z)=\frac{1}{2}\sin{(2z+B)}+A,$$ where $A$, $B$ are constants and $c=k\pi+\frac{\pi}{2}$ $(k\in \mathds{Z})$.
\par (3) If $q=-1$, then Eq.\eqref{E:2.8} has the finite order of transcendental entire solution $$w(z)=\frac{1}{2}\sin(2z+B)+A,$$  where $A$, $B$ are constants  and $B+c=k\pi$ $(k\in \mathds{Z})$.
\end{theorem}
\begin{rem}
This theorem implies that if Eq.\eqref{E:2.8} has the finite order of transcendental entire solutions, then Eq.\eqref{E:2.8} reduce to Eq.\eqref{E:2.4}, or reduce to the following equation
\begin{align*}
w^{\prime}(z)^{2}+[w(-z+c)-w(z)]^{2}=1.
\end{align*}
\end{rem}

\begin{theorem}  Let $m$, $n$ be positive integers. If $m>n=1$  or $n>m=1$, then  Eq.\eqref{E:2.5} does not posses any zero-order of transcendental entire solution.
\end{theorem}
When $m=n=1$ and $q=1$ in Eq.\eqref{E:2.5}, then $w^{\prime}(z)+w(z+c)=1$ can admit a transcendental entire solution (see \cite{s9}).

\begin{theorem}  Let $m$, $n$ be positive integers. If $n>m=1$, then Eq.\eqref{E:2.6} does not posses any zero-order of transcendental entire solution.
\end{theorem}

\section{Results on the Fermat type of systems of complex differential-difference equations}
\par In 2016, Gao \cite{s11} considered the existence of transcendental entire solutions of the Fermat type of systems of  complex differential-difference equations as bellow:
\begin{align}\label{E:3.1}  \left\{
\begin{array}{ll}w^{\prime}_{1}(z)^{2}+w_{2}(z+c)^{2}=1,\\
w^{\prime}_{2}(z)^{2}+w_{1}(z+c)^{2}=1.~~~~
         \end{array}\right.
\end{align}
{\bf Theorem F. (\cite{s11})} Let $(w_{1}(z),w_{2}(z))$ be a finite order of transcendental entire solutions of Eq.\eqref{E:3.1}, then $(w_{1}(z),w_{2}(z))$ must satisfy $$(w_{1}(z),w_{2}(z))=(\sin(z+B_{1}),\ \sin{(z+B_{2})}),$$ where $B_{1}$ and $B_{2}$ are constants, $c=k\pi$, $B_{1}-B_{2}=l\pi$, $l$, $k$ are integers.

 However, this result is not complete, which complete version will be included in Theorem 3.2 later.
\par Gao \cite{s11} also researched the more general case as follow.
\begin{align}\label{E:3.01}  \left\{
\begin{array}{ll}w^{\prime}_{1}(z)^{n_{1}}+w_{2}(z+c)^{m_{1}}=1,\\
w^{\prime}_{2}(z)^{n_{2}}+w_{1}(z+c)^{m_{2}}=1.
         \end{array}\right.
\end{align}
where $c$ is a complex constant, $n_{1}, n_{2}, m_{1}, m_{2}$ are positive integers and $n_{i}>1,~i=1,2$. He proved the following results.\\
{\bf Theorem G. (\cite{s11})} The Eq.\eqref{E:3.01} has no transcendental entire solution with finite order, if it satisfies  one of the following two conditions:
\par (i) $m_{1}m_{2}>n_{1}n_{2}$;
\par(ii)$m_{i}>\frac{n_{i}}{n_{i}-1}, ~~i=1,2.$\\
\par Enlightened by \cite{s11}, we will mainly investigate the following two types of systems of complex differential-difference equations.
\begin{align}\label{E:3.2} \left\{
\begin{array}{ll}w^{\prime}_{1}(z)^{n_{1}}+w_{2}(qz+c)^{m_{1}}=1,\\
w^{\prime}_{2}(z)^{n_{2}}+w_{1}(qz+c)^{m_{2}}=1,~~~~~~~~~~~~
         \end{array}\right.
\end{align}
where $q(\neq0)$ and $c$ are complex constants.
\par The order of solutions $(w_{1}(z),w_{2}(z))$ of Eq.\eqref{E:3.2}, denote
$$\rho = \rho(w_{1},w_{2})=\max\{\rho(w_{1}),\rho(w_{2})\},$$
where
$$\rho(w_{k})=\limsup_{r\rightarrow \infty}\frac{\log^{+}T(r,w_{k})}{\log{r}},~~~k=1,2.$$
\par We state our result as  follows
\begin{theorem}  Let $n_{1}$, $m_{1}$, $n_{2}$, $m_{2}$ be positive integers satisfying $n_{1}\geq2$ $m_{1}\geq2$, and at lest one is not equal to 2, or $n_{2}\geq2$, $m_{2}\geq2$ and at lest one is not equal to 2 . Then there are no non-constant entire solutions $(w_{1}(z),w_{2}(z))$ that satisfy Eq.\eqref{E:3.2}.
\end{theorem}
When $n_{1}=m_{1}=n_{2}=m_{2}=2$ , Eq.\eqref{E:3.2} take the form:
\begin{align}\label{E:3.4} \left\{
\begin{array}{ll}w^{\prime}_{1}(z)^{2}+w_{2}(qz+c)^{2}=1,\\
w^{\prime}_{2}(z)^{2}+w_{1}(qz+c)^{2}=1,~~~~~~~~~~~~
         \end{array}\right.
\end{align}
where $q(\neq0)$ and $c$ are complex constants.
\par Considering the different values $q$, we divide into three case. Precisely,
\begin{theorem} (1) If $q\neq \pm 1$, then Eq.\eqref{E:3.4} does not posses any  finite order of transcendental entire solutions.
\par (2) If $q=1$, then Eq.\eqref{E:3.4} has the finite order of transcendental entire solutions $(w_{1}(z),w_{2}(z))$ with $$(w_{1}(z),w_{2}(z))=(\sin{(z+B_{1})}, \ \sin{(z+B_{2})}),$$
where  $B_{1}$, $B_{2}$ are two complex constants, $c=k\pi$ and $B_{1}-B_{2}=l\pi$, or $c=\frac{\pi}{2}+k\pi$ and $B_{1}-B_{2}=\frac{\pi}{2}+l\pi$  $(l, k\in \mathds{Z})$.
\par (3) If $q=-1$, then Eq.\eqref{E:3.4} has the finite order of transcendental entire solutions $(w_{1}(z),w_{2}(z))$ with $$(w_{1}(z),w_{2}(z))=(\sin{(z+B_{1})},\ \sin{(z+B_{2})}),$$
where  $B_{1}$, $B_{2}$ are two complex constants, $B_{1}+B_{2}+c=k\pi$ $(k\in \mathds{Z})$.
\end{theorem}
\begin{rem}Theorem 3.2 is show that if Eq.\eqref{E:3.4} has the finite order of transcendental entire solutions, then Eq.\eqref{E:3.4} maybe reduce to Eq.\eqref{E:3.1} or the following systems of complex differential-difference equation
\begin{align}\nonumber \left\{
\begin{array}{ll}w^{\prime}_{1}(z)^{2}+w_{2}(-z+c)^{2}=1,\\
w^{\prime}_{2}(z)^{2}+w_{1}(-z+c)^{2}=1.
         \end{array}\right.
\end{align}
\end{rem}

\section{Lemmas for the proof of the theorems}
In this section, we introduce some lemmas which ere essential in the process of proving our results.
\begin{lemma} (see\cite{s2}) Let $w(z)$ be a finite order $\rho(w)$  of transcendental entire function with a zero of multiplicity $k\geq 0$ at $z=0$. Let the other zeros of $w(z)$ be at  $z_{1},z_{2}\cdots,$ each zero being repeated as many as its multiplicity implies.  then
\begin{align*}
w(z)=z^{k}P(z)e^{Q(z)},
\end{align*}
where $P(z)$ is the canonical product for all non-zero zero points of $w(z)$, $Q(z)$ is a polynomial and $\deg{Q(z)}\leq \rho(w)$.
\end{lemma}
Lemma 4.1 is called the Hadamard factorization theorem of entire function. the following result is a key lemma in this paper.
\begin{lemma}(see\cite{s2})  Let
$w_{j}(z)$ be meromorphic functions, $w_{k}(z)(k=1,2,\cdots, n-1)$  be not constants, satisfy $\sum_{j=1}^{n}w_{j}=1$ and $n\geq3$, If $w_{n}(z)\not\equiv 0$ and
\begin{align*}
\sum^{n}_{j=1}N(r,\frac{1}{w_{j}})+(n-1)\sum^{n}_{j=1}\overline{N}(r,w_{j})<(\lambda+o(1))T(r,w_{k}),
\end{align*}
where $\lambda<1$, $k=1,2,\cdots,n-1$, then $w_{n}\equiv 1$.
\end{lemma}
\begin{lemma}(see\cite{s13}) Let $w(z)$ be a non-constant zero-order meromorphic function, and
$q \in C\setminus \{0\}$. Then
\begin{align*}
m\bigg(r,\frac{w(qz)}{w(z)}\bigg)=o(T(r,w))~~(r \rightarrow \infty, r \not\in E)
\end{align*}
where $E$ is a set of logarithmic density 0.
\end{lemma}
\begin{lemma}(see\cite{s14})Let $w(z)$ be a transcendental meromorphic solutions of
 \begin{align*}
f^{n}P(z,w)=Q(z,w),
\end{align*}
where $P(z,w)$ and $Q(z,w)$ are polynomial in $w$ and its derivatives with meromorphic coefficients, say $\{ a_{\lambda}\mid \lambda\in I\}$, such that $m(r,a_{\lambda})=S(r,w)$ for all $\lambda\in I$. If the total degree of $Q(z,w)$ as a polynomial in $w$ and its derivatives is $\leq n$, then
$$m(r,P(z,w))=S(r,w).$$
\end{lemma}
\begin{lemma}(see\cite{s8})
Let $w(z)$ be a meromorphic function of finite
order $\rho(w)$. If $\varepsilon>0$, then
 \begin{align*}
T(r,w(z+c))=T(r,w)+O(r^{\rho(w)-1+\varepsilon})+O(\log{r}).
\end{align*}
\end{lemma}

\section{The proof of Theorems 2.1-2.5}
\emph{\bf Proof of Theorem 2.1:}
 According to \eqref{E:2.7}, let
$$f(z)=w^{\prime}(z),~~~~~~g(z)=w(qz+c),$$
then
\begin{align*}
f(z)^{n}+g(z)^{m}=1.
\end{align*}
Because $n\geq2$, $m\geq2$  and at lest one is not equal to 2, we have
$$\frac{1}{n}+\frac{1}{m}<1.$$
By Theorem A, we can conclude that Eq.\eqref{E:2.5} does not posses non-constant entire solutions.
\par Similarly, Eq.\eqref{E:2.6} does not posses non-constant entire solutions.$\hfill\Box$\\

\emph{\bf Proof of Theorem 2.2:} Let $w(z)$ be a transcendental entire solution with finite order of
 Eq.\eqref{E:2.7}, then
 \begin{align*} (w^{\prime}(z)+iw(qz+c))(w^{\prime}(z)-iw(qz+c))=1.
\end{align*}
Thus, $w^{\prime}(z)+iw(qz+c)$ and $w^{\prime}(z)-iw(qz+c)$ have no zeros. Applying Lemma 4.1, we must have
\begin{align} \label{E:5.1}\left\{
\begin{array}{ll}w^{\prime}(z)+iw(qz+c)=e^{p(z)},\\
w^{\prime}(z)-iw(qz+c)=e^{-p(z)},
         \end{array}\right.
\end{align}
where $p(z)$ is a nonconstant polynomial.  By \eqref{E:5.1},  we gain
\begin{align}
\label{E:5.2}&w^{\prime}(z)=\frac{e^{p(z)}+e^{-p(z)}}{2},
\end{align}
and
\begin{align}
\label{E:5.3}&w(qz+c)=\frac{e^{p(z)}-e^{-p(z)}}{2i}.
\end{align}
Combining \eqref{E:5.2} with \eqref{E:5.3}, we get
\begin{align}
\label{E:5.4}&w^{\prime}(qz+c)=\frac{e^{p(qz+c)}+e^{-p(qz+c)}}{2}=\frac{p^{\prime}(z)e^{p(z)}+p^{\prime}(z)e^{-p(z)}}{2i}.
\end{align}
By \eqref{E:5.4}, we gain
\begin{align}\label{E:5.5}
\frac{qi}{p^{\prime}(z)}\big(e^{p(qz+c)+p(z)}+e^{p(z)-p(qz+c)}\big)-e^{2p(z)}=1.
\end{align}
Next we will confirm that $e^{2p(z)}$ is not a constant. Otherwise, suppose $e^{2p(z)} \equiv d_{1}$, where $d_{1}$ is a constant, then we have $p(z)$ is a non-zero constant, which shows that $w(z)$ must be a constant, a contradiction. Meanwhile, we will also affirm that $\frac{qi}{p^{\prime}(z)}e^{p(qz+c)+p(z)}$ and $\frac{qi}{p^{\prime}(z)}e^{p(z)-p(qz+c)}$  cannot be all constants. Otherwise, assume $\frac{qi}{p^{\prime}(z)}e^{p(qz+c)+p(z)}=d_{2}$, where $d_{2}$ is a constant. Moreover, it implies that ${p^{\prime}(z)}$ has no zeros. Let $p(z)=az+b$, where $a$ is a nonzero constant, $b$ is a constant. Because $\frac{qi}{p^{\prime}(z)}e^{p(qz+c)+p(z)}=d_{2}$, thus $p(qz+c)+p(z)$ is a constant, that is $$p(qz+c)+p(z)=a(qz+c)+b+az+b=(qa+a)z+2b+ac$$ is a constant.  Then $q=-1$, and we get \begin{align*}\frac{qi}{p^{\prime}(z)}e^{p(z)-p(qz+c)}&=\frac{-i}{a}e^{p(z)-p(qz+c)}=\frac{-i}{a}e^{az+b-(a(-z+c)+b)}\\
&=\frac{-i}{a}e^{2az-ac}.
\end{align*}
Obviously, $\frac{qi}{p^{\prime}(z)}e^{p(z)-p(qz+c)}$ is not a constant. Similarly, we can prove that $\frac{qi}{p^{\prime}(z)}e^{p(qz+c)+p(z)}$ is not a constant, if $\frac{qi}{p^{\prime}(z)}e^{p(z)-p(qz+c)}$ be a constant. In the following, we will discuss two cases.
\par\emph{Case (i):} Let $e^{2p(z)}$ and $\frac{qi}{p^{\prime}(z)}e^{p(z)+p(qz+c)}$  be not constants. Applying Lemma 4.2 and \eqref{E:5.5}, we get that
\begin{align}\label{E:5.6}\frac{qi}{p^{\prime}(z)}e^{p(z)-p(qz+c)}\equiv 1,\end{align} hence ${p^{\prime}(z)}$ has no zeros and $p(z)$ is a nonconstant polynomial, then $\deg{p(z)}=1$. Assume $p(z)=az+b$, where $a$ is a nonzero constant, $b$ is a constant. By \eqref{E:5.6}, we obtain $a-aq=0$, that is $q=1$. Moreover, the Fermat type equation \eqref{E:2.7} reduce to \eqref{E:2.3}.
\par \emph{Case(ii):} Let $e^{2p(z)}$ and $\frac{qi}{p^{\prime}(z)}e^{p(z)-p(qz+c)}$  be not constants. Applying Lemma 4.2 and \eqref{E:5.5}, we get that
\begin{align}\label{E:5.7}\frac{qi}{p^{\prime}(z)}e^{p(z)+p(qz+c)}\equiv 1.\end{align}
 Hence ${p^{\prime}(z)}$ has no zeros and $p(z)$ is a nonconstant polynomial, then $\deg{p(z)}=1$. Assume $p(z)=az+b$, where $a$ is a nonzero constant, $b$ is a constant. By \eqref{E:5.7}, we obtain $a+aq=0$, that is $q=-1$. Thus, we obtain
 \begin{align}\label{E:5.8}\frac{qi}{p^{\prime}(z)}e^{p(z)+p(qz+c)}=\frac{-i}{a}e^{2b+ac}= 1.
\end{align}
Hence $e^{2b+ac}=ai$. From \eqref{E:5.2}, we gain
\begin{align}\label{E:5.9}
w(z)=\frac{\frac{1}{a}e^{az+b}-\frac{1}{a}e^{-az-b}}{2}+A,
\end{align}
where $A$ is a constant.
By \eqref{E:5.3} and \eqref{E:5.8}, we get
\begin{align*}
w(z)&=\frac{e^{a(-z+c)+b}-e^{-a(-z+c)-b}}{2i}=\frac{e^{-az-b+ac+2b}-e^{az+b-ac-2b}}{2i}.
\end{align*}
That is
\begin{align}
\label{E:5.10}w(z)=\frac{ae^{-az-b}+\frac{1}{a}e^{az+b}}{2}.
\end{align}
Comparing \eqref{E:5.9} with \eqref{E:5.10}, we obtain $a=\pm i$ and $A=0$.
\par\emph{Subcase(i):} If $a=i$, then $2b+ic=(2k+1)\pi i,$ that is $2bi-c=(2k+1)\pi$. From \eqref{E:5.10}, we get
\begin{align}\nonumber
w(z)&=\frac{ae^{-az-b}+\frac{1}{a}e^{az+b}}{2}=\frac{ie^{-iz-b}+\frac{1}{i}e^{iz+b}}{2}\\
\nonumber&=\frac{e^{i(z-bi)}-e^{-i(z-bi)}}{2i}\\
\label{E:5.11}&=\sin(z-ib).
\end{align}
Let $B=-bi$, then $w(z)=\sin(z+B)$, where $B$ is also a constant, and $2B+c=(2k+1)\pi$.
\par \emph{Subcase(ii):} If $a=-i$, then $2b-ic=2k\pi i,$ that is $2bi+c=2k\pi$. From \eqref{E:5.10}, we obtain
\begin{align}\nonumber
w(z)&=\frac{ae^{-az-b}+\frac{1}{a}e^{az+b}}{2}=\frac{-ie^{iz-b}+\frac{1}{-i}e^{-iz+b}}{2}\\
\nonumber&=\frac{e^{i(z+ib)}-e^{-i(z+ib)}}{2i}\\\label{E:5.12}&=\sin(z+ib).
\end{align}
Let $B=bi$, then $w(z)=\sin(z+B)$, where $B$ is also a constant, and $2B+c=2k\pi$. $\hfill\Box$

\emph{\bf Proof of Theorem 2.3:} Let $w(z)$ be a transcendental entire solution with finite order of
 the complex differential-difference equations \eqref{E:2.8}. Let $g(z)=w(qz+c)-w(z)$, then
 \begin{align*} (w^{\prime}(z)+ig(z))(w^{\prime}(z)-ig(z))=1.
\end{align*}
Thus, $w^{\prime}(z)+ig(z)$ and $w^{\prime}(z)-ig(z)$ have no zeros. Applying Lemma 4.1, we can assume that
\begin{align} \label{E:5.13}\left\{
\begin{array}{ll}w^{\prime}(z)+ig(z)=e^{p(z)},\\
w^{\prime}(z)-ig(z)=e^{-p(z)},
         \end{array}\right.
\end{align}
where $p(z)$ is a nonconstant polynomial.  From \eqref{E:5.13},  we get
\begin{align}
\label{E:5.14}&w^{\prime}(z)=\frac{e^{p(z)}+e^{-p(z)}}{2},
\end{align}
and
\begin{align}
\label{E:5.15}g(z)=\frac{e^{p(z)}-e^{-p(z)}}{2i}.
\end{align}
From \eqref{E:5.14} and \eqref{E:5.15}, we obtain
\begin{align}
\label{E:5.16}g^{\prime}(z)=qw^{\prime}(qz+c)-w^{\prime}(z)=\frac{p^{\prime}(z)(e^{p(z)}+e^{-p(z)})}{2i}.
\end{align}
Then
\begin{align*}
w^{\prime}(qz+c)&=\frac{p^{\prime}(z)(e^{p(z)}+e^{-p(z)})}{2qi}+\frac{e^{p(z)}+e^{-p(z)}}{2q}\\
&=\frac{(p^{\prime}(z)+i)(e^{p(z)}+e^{-p(z)})}{2qi}.
\end{align*}
That is
\begin{align}
\label{E:5.17}w^{\prime}(qz+c)=\frac{(p^{\prime}(z)+i)(e^{p(z)}+e^{-p(z)})}{2qi}.
\end{align}
By \eqref{E:5.14}, we get
\begin{align}
\label{E:5.18}w^{\prime}(qz+c)=\frac{e^{p(qz+c)}+e^{-p(qz+c)}}{2}.\end{align}
Therefore, from \eqref{E:5.17} and \eqref{E:5.18}, we gain
\begin{align}\label{E:5.19}
\frac{p^{\prime}(z)+i}{qi}\big(e^{p(qz+c)+p(z)}+e^{p(qz+c)-p(z)}\big)-e^{2p(qz+c)}=1.
\end{align}

Similar in the proof of Theorem 2.2, we confirm that $e^{2p(qz+c)}$ is not constant, and
$\frac{p^{\prime}(z)+i}{qi}$ $e^{p(qz+c)-p(z)}$ and
$\frac{p^{\prime}(z)+i}{qi}e^{p(qz+c)+p(z)}$ are not constants simultaneously.
Thus, we need discuss the following two cases.
\par \emph{Case (i):} Let $e^{2p(qz+c)}$ and $\frac{p^{\prime}(z)+i}{qi}e^{p(qz+c)+p(z)}$  be not constants. Applying Lemma 4.2 and \eqref{E:5.5} , we get
\begin{align}\label{E:5.20}
\frac{p^{\prime}(z)+i}{qi}e^{p(qz+c)-p(z)}\equiv 1,
\end{align}
hence $\frac{p^{\prime}(z)+i}{qi}$ has no zeros and $p(z)$ is a nonconstant polynomial. Moreover $\deg{p(z)}=1$. Assume $p(z)=az+b$, where $a$ is a nonzero constant, $b$ is a constant. By \eqref{E:5.20}, we obtain $a-aq=0$, that is $q=1$. Moreover, the Fermat type equation \eqref{E:2.8} reduces to \eqref{E:2.4}.

\par\emph{Case(ii):} Let $e^{2p(qz+c)}$ and $\frac{p^{\prime}(z)+i}{qi}e^{p(qz+c)-p(z)}$  be not constants. Applying Lemma 4.2 and \eqref{E:5.19}, we have
\begin{align}\label{E:5.21}
\frac{p^{\prime}(z)+i}{qi}e^{p(z)+p(qz+c)}\equiv 1,
\end{align}
 hence $\frac{p^{\prime}(z)+i}{qi}$ has no zeros and $p(z)$ is a nonconstant polynomial, then $\deg{p(z)}=1$. Assume $p(z)=az+b$, where $a$ is a nonzero constant, $b$ is a constant. By \eqref{E:5.21}, we obtain $a+aq=0$, that is $q=-1$. Thus we derive
 \begin{align}\label{E:5.22}\frac{p^{\prime}(z)+i}{qi}e^{p(z)+p(qz+c)}=\frac{a+i}{-i}e^{2b+ac}= 1.
\end{align}
That is  $e^{2b+ac}=\frac {1}{ai-1}$. From \eqref{E:5.15}, we gain
\begin{align}\label{E:5.23}
w(z)=\frac{\frac{1}{a}e^{az+b}-\frac{1}{a}e^{-az-b}}{2}+A,
\end{align}
where $A$ is a constant.
By \eqref{E:5.14}, we get
\begin{align*}
w(-z+c)&=\frac{e^{az+b}-e^{-az-b}}{2i}+w(z)\\
&=\frac{(\frac{i}{a}+1)(e^{az+b}-e^{-az-b})}{2i}+A.
\end{align*}
That is
\begin{align}
\nonumber w(z)&=\frac{(\frac{i}{a}+1)(e^{a(-z+c)+b}-e^{-a(-z+c)-b})}{2i}+A\\
\nonumber&=\frac{(\frac{i}{a}+1)(e^{-az-b+(ac+2b)}-e^{az+b-(ac+2b)})}{2i}+A
\end{align}
\begin{align}
\nonumber&=\frac{\frac{i+a}{a}\frac {1}{ai-1}e^{-az-b}-\frac{i+a}{a}(ai-1)e^{az+b}}{2i}+A\\
\nonumber&=\frac{\frac{1}{ai}e^{-az-b}-\frac{i(i+a)^{2}}{a}e^{az+b}}{2i}+A\\
\label{E:5.24}&=\frac{-\frac{1}{a}e^{-az-b}-\frac{(i+a)^{2}}{a}e^{az+b}}{2}+A.
\end{align}
Comparing \eqref{E:5.23} with \eqref{E:5.24}, we obtain $-\frac{(i+a)^{2}}{a}=\frac{1}{a}$. Then $a=0$ or $a=-2i$. Obviously,the case $a=0$ should be desert, which is contradiction with our hypothesis. Therefore, we only consider that $a=-2i$ and $e^{2b-2ic}=1$. From \eqref{E:5.23}, we have
\begin{align*}
w(z)&=\frac{1}{2}\bigg(\frac{e^{2iz-b}-e^{-2iz+b}}{2i}\bigg)+A\\
&=\frac{1}{2}\sin{(2z+ib)}+A.
\end{align*}
Let $B=bi$, then $w(z)=\frac{1}{2}\sin(2z+B)+A$, where $A$, $B$ are constants, and $B+c=k\pi$.$\hfill\Box$

\emph{\bf Proof of Theorem 2.4:} Let $w(z)$ be a zero order of transcendental entire solution  of
 Eq.\eqref{E:2.5}.

(i) If $n>m=1$, that is
\begin{align}\label{E:5.26} w^{\prime}(z)^{n}+w(qz+c)=1.
\end{align}
By differentiating \eqref{E:5.26}, we gain
\begin{align} \label{E:5.27}
nw^{\prime}(z)^{n-1}w^{\prime\prime}=-qw^{\prime}(qz+c).
\end{align}
Let $F(z)=w^{\prime}(z)$, thus we get by \eqref{E:5.27}
\begin{align} \label{E:5.28}
nF(z)^{n-1}F^{\prime}(z)=-a(z)F(z),
\end{align}
where $a(z)=\frac{q F(qz+c)}{F(z)}$ is a meromorphic function of $z$. By Lemma 4.3 and \eqref{E:5.28}, we get $m(r,a(z))=S(r,w^{\prime}(z))$.
Applying Lemma 4.4, we have
$$m(r,F^{\prime}(z))=m(r,w^{\prime\prime}(z))=S(r,w^{\prime}(z)),$$
which contradicts the fact that $w(z)$ being a zero order of transcendental entire function.

(ii) If $m>n=1$, that is
\begin{align}\label{E:5.29} w^{\prime}(z)+w(qz+c)^{m}=1.
\end{align}
From \eqref{E:5.29} and Lemma 4.5, we get
\begin{align} \nonumber
mT(r,w)+S(r,w)&=mT(r,w(qz+c))\\
\nonumber&=T(r,w^{\prime}-1)+S(r,w)\\
\nonumber &=T(r,w^{\prime})+S(r,w)\\
\nonumber&\leq T(r,w)+S(r,w).
\end{align}
That is
 \begin{align} \nonumber
(m-1)T(r,w)\leq S(r,w).
\end{align}
Thus, in view of $m>1$, we gain  $T(r,w)=S(r,w)$, a contradiction. So our theorem follows. $\hfill\Box$

\emph{\bf Proof of Theorem 2.5:} Let $w(z)$ be a zero order of transcendental entire solution  of
 the complex differential-difference equations \eqref{E:2.6}.
If $n>m=1$, that is
\begin{align}\label{E:5.30} w^{\prime}(z)^{n}+w(qz+c)-w(z)=1.
\end{align}
Differentiating \eqref{E:5.30}, we gain
\begin{align} \label{E:5.31}
nw^{\prime}(z)^{n-1}w^{\prime\prime}=-qw^{\prime}(qz+c)+w^{\prime}(z).
\end{align}
Let $F(z)=w^{\prime}(z)$, and \eqref{E:5.31}, thus we get
\begin{align} \label{E:5.32}
nF(z)^{n-1}F^{\prime}(z)=(a(z)+1)F(z),
\end{align}
where $a(z)=\frac{q F(qz+c)}{F(z)}$ is a meromorphic function of $z$. By Lemma 4.3 and \eqref{E:5.32}, we get $m(r,a(z))=S(r,w^{\prime}(z))$.
Applying Lemma 4.4, we gain
$$m(r,F^{\prime}(z))=m(r,w^{\prime\prime}(z))=S(r,w^{\prime}(z)).$$
It is a contradiction with the fact $w(z)$ being a zero order of transcendental entire function. $\hfill\Box$

\section{The proof of Theorem 3.1 and Theorem 3.2}
\emph{\bf Proof of Theorem 3.1:} If $m_{1}\geq2$, $n_{1}\geq 2 $ and at lest one is not equal to 2, in view of the first equation of Eq.\eqref{E:3.2}, assuming
$$f(z)=w_{1}^{\prime}(z),~~~~~~g(z)=w_{2}(qz+c),$$  we have
\begin{align*}
f(z)^{n_{1}}+g(z)^{m_{1}}=1,
\end{align*}
and
\begin{align*}
\frac{1}{n_{1}}+\frac{1}{m_{1}}<1.
\end{align*}
From Theorem A, thus the first equation of the systems of complex differential-difference equations has no entire solutions.
If $m_{2}\geq2$, $n_{2}\geq 2 $ and at lest one is not equal to 2, then the second equation of Eq.\eqref{E:3.2}has no entire solutions.  Therefore, Eq.\eqref{E:3.2} does not posses any nonconstant entire solutions.$\hfill\Box$\\

\par\emph{\bf Proof of Theorem 3.2:} Assume $(w_{1}(z),w_{2}(z))$ be a finite order of transcendental entire solution of the Fermat type of the systems of complex differential-difference equations  \eqref{E:3.4}. From \eqref{E:3.4}, we get
\begin{align}\label{E:6.1}\left\{
\begin{array}{ll}(w^{\prime}_{1}(z)+iw_{2}(qz+c))(w^{\prime}_{1}(z)-iw_{2}(qz+c))=1,\\
(w^{\prime}_{2}(z)+iw_{1}(qz+c))(w^{\prime}_{2}(z)-iw_{1}(qz+c))=1.
         \end{array}\right.
\end{align}
Thus, $w^{\prime}_{1}(z)+iw_{2}(qz+c)$, $w^{\prime}_{1}(z)-iw_{2}(qz+c)$, $w^{\prime}_{2}(z)+iw_{1}(qz+c)$ and $w^{\prime}_{2}(z)-iw_{1}(qz+c)$ have no zeros. Combining  \eqref{E:6.1} and Lemma 4.1, we can assume that
\begin{align}\label{E:6.2} \left\{
\begin{array}{ll}w^{\prime}_{1}(z)+iw_{2}(qz+c)=e^{p(z)},\\
w^{\prime}_{1}(z)-iw_{2}(qz+c)=e^{-p(z)},\\
w^{\prime}_{2}(z)+iw_{1}(qz+c)=e^{h(z)},\\
w^{\prime}_{2}(z)-iw_{1}(qz+c)=e^{-h(z)},
         \end{array}\right.
\end{align}
where $p(z)$ and $h(z)$ are nonconstant polynomials.
From \eqref{E:6.2}, we get
\begin{align}
\label{E:6.3}&w^{\prime}_{1}(z)=\frac{e^{p(z)}+e^{-p(z)}}{2},\\
\label{E:6.4}&w_{2}(qz+c)=\frac{e^{p(z)}-e^{-p(z)}}{2i},\\
\label{E:6.5}&w^{\prime}_{2}(z)=\frac{e^{h(z)}+e^{-h(z)}}{2},\\
\label{E:6.6}&w_{1}(qz+c)=\frac{e^{h(z)}-e^{-h(z)}}{2i}.
\end{align}
Combining \eqref{E:6.3} and \eqref{E:6.6}, \eqref{E:6.4} and \eqref{E:6.5} respectively, we obtain
\begin{align}\label{E:6.7} \left\{
\begin{array}{ll}w^{\prime}_{1}(z)=\frac{e^{p(z)}+e^{-p(z)}}{2},\\
w_{1}(qz+c)=\frac{e^{h(z)}-e^{-h(z)}}{2i}.
         \end{array}\right.
\end{align}
and
\begin{align}\label{E:6.8}  \left\{
\begin{array}{ll}w^{\prime}_{2}(z)=\frac{e^{h(z)}+e^{-h(z)}}{2},\\
w_{2}(qz+c)=\frac{e^{p(z)}-e^{-p(z)}}{2i}.
         \end{array}\right.
\end{align}
From \eqref{E:6.7} and \eqref{E:6.8}, we get
\begin{align}\label{E:6.9}
\frac{qi}{h^{\prime}(z)}(e^{p(qz+c)+h(z)}+e^{h(z)-p(qz+c)})-e^{2h(z)}=1,
\end{align}
and
\begin{align}\label{E:6.10}
\frac{qi}{p^{\prime}(z)} (e^{h(qz+c)+p(z)}+e^{p(z)-h(qz+c)})-e^{2p(z)}=1.
\end{align}

\par We claim that $e^{2p(z)}$ and $e^{2h(z)}$ are not constants,  and $\frac{q i}{h^{\prime}(z)}e^{h(z)+p(qz+c)}$ and $\frac{q i}{h^{\prime}(z)}e^{h(z)-p(qz+c)}$, $\frac{q i}{p^{\prime}(z)}e^{p(z)+h(qz+c)}$ and $\frac{q i}{p^{\prime}(z)}e^{p(z)-h(qz+c)}$ cannot be all constants respectively.
\par\emph{Case (i):} If $e^{2p(z)}\equiv d_{5}$ and $e^{2h(z)}\equiv d_{6}$, where $d_{5}$ and $d_{6}$ are complex constants, then $p(z)$ and $h(z)$ are constants. Obviously, $w_{1}(z)$ and $w_{2}(z)$ are constants, which are contradiction.
\par\emph{Case (ii):} If $e^{2p(z)}\equiv d_{7}$ and $e^{2h(z)}\not\equiv d_{8}$, where $d_{7}$ and $d_{8}$ are complex constants, then $p(z)$ are constants. Thus, by \eqref{E:6.3}, we get  $w_{1}(z)=az+b_{1}$, where $a$ is a nonzero constant. It is a contradiction with the second of the systems of equation \eqref{E:6.7}.
\par\emph{Case (iii):} If $e^{2p(z)}\not\equiv d_{9}$ and $e^{2h(z)}\equiv d_{10}$, where $d_{9}$ and $d_{10}$ are complex constants, we can also get a contradiction with the second of the systems of equation \eqref{E:6.8}.
\par Hence, $e^{2p(z)}$ and $e^{2h(z)}$ are not constants.
\par\emph{Case (iv):} We will affirm that $\frac{q i}{h^{\prime}(z)}e^{h(z)+p(qz+c)}$ and $\frac{q i}{h^{\prime}(z)}e^{h(z)-p(qz+c)}$ cannot be all constants simultaneously. If $\frac{q i}{h^{\prime}(z)}e^{h(z)+p(qz+c)}\equiv d_{11}$ and $\frac{q i}{h^{\prime}(z)}e^{h(z)-p(qz+c)}\equiv d_{12}$, where $d_{11}$ and $d_{12}$ are complex constants. Thus, $h^{\prime}(z)$ has no zeros, then for any constant $b_{2}$, $h(z)=az+b_{2}$, where $a$ is a nonzero constant.  By $\frac{q i}{h^{\prime}(z)}e^{h(z)-p(qz+c)}\equiv d_{12}$, we assume that $p(z)=\frac{a}{q}z+b_{1}$ and $h(z)=az+b_{2}$, where $b_{1}$ is any constant. Then
$$\frac{q i}{h^{\prime}(z)}e^{h(z)+p(qz+c)}=\frac{q i}{a}e^{az+b_{2}+\frac{a}{q}(qz+c)+b_{1}}=\frac{q i}{a}e^{2az+b_{2}+\frac{ac}{q}+b_{1}},$$
which contradicts with $\frac{q i}{h^{\prime}(z)}e^{h(z)+p(qz+c)}\equiv d_{11}$. Thus, $\frac{q i}{h^{\prime}(z)}e^{h(z)+p(qz+c)}$ and $\frac{q i}{h^{\prime}(z)}e^{h(z)-p(qz+c)}$ cannot be constants simultaneously.
\par\emph{Case (v):} Similarly, we also can confirm that $\frac{q i}{p^{\prime}(z)}e^{p(z)+h(qz+c)}$ and $\frac{q i}{p^{\prime}(z)}e^{p(z)-h(qz+c)}$ cannot be constants simultaneously.\\

\par According to the above discussion, the following four cases, sixteen subcases need to be investigated.

\par {\bf Case (A):} Let $e^{2p(z)}$, $e^{2h(z)}$, $\frac{q i}{h^{\prime}(z)}e^{h(z)+p(qz+c)}$ and $\frac{q i}{p^{\prime}(z)}e^{p(z)+h(qz+c)}$  be not constants. From \eqref{E:6.9}, \eqref{E:6.10} and Lemma 4.2, we get
\begin{align}\label{E:6.11}
\frac{q i}{h^{\prime}(z)}e^{h(z)-p(qz+c)}\equiv 1,
\end{align}
and
\begin{align}\label{E:6.12}
\frac{q i}{p^{\prime}(z)}e^{p(z)-h(qz+c)}\equiv 1.
\end{align}
Thus, $p^{\prime}(z)$ and $h^{\prime}(z)$ have no zeros, $p(z)$ and $h(z)$ are nonconstant polynomials, which implies that $\deg{p(z)}=1$ and $\deg{h(z)}=1$.
Hence, let $p(z)=a_{1}z+b_{1}$ and $h(z)=a_{2}z+b_{2}$, where $a_{1}(\neq0)$, $a_{2}(\neq0)$, $b_{1}$ and $b_{2}$ are constants. From \eqref{E:6.11}, we get
$$\frac{q i}{h^{\prime}(z)}e^{h(z)-p(qz+c)}=\frac{q i}{a_{2}}e^{a_{2}z+b_{2}-a_{1}(qz+c)+b_{1}}\equiv 1,$$
Thus, we can deduce that $a_{2}=qa_{1}$.
\par Similarly, from \eqref{E:6.12}, we can gain $a_{1}=qa_{2}$.  Then $a_{1}a_{2}=q^{2}a_{1}a_{2}$, that is  $q=1$ or $q=-1$.

\par Subcase (i): $q=1$. For the case $q=1$ in the systems of equations \eqref{E:3.4}, Gao has already discussed it in \cite{s11}, however, his investigation is not complete, which not including Subcase (i) in Case (B) in the later part of this paper. For the completeness of the proof, here we also show its proving process.

   If $q=1$, then $a_{1}=a_{2}$, so we set $p(z)=az+b_{1}$ and $h(z)=az+b_{2}$. From \eqref{E:6.11} and \eqref{E:6.12} we have
\begin{align}
e^{b_{2}-ac-b_{1}}=-ia, \label{E:6.12a}
\end{align}
\begin{align}
e^{b_{1}-ac-b_{2}}=-ia. \label{E:6.12b}
\end{align}
By \eqref{E:6.6}, then
\begin{align*}
w_{1}(z+c)=\frac{e^{az+b_{2}}-e^{-az-b_{2}}}{2i}.
\end{align*}
Combining with \eqref{E:6.12a}, thus
\begin{align*}
w_{1}(z)&=\frac{e^{a(z-c)+b_{2}}-e^{-a(z-c)-b_{2}}}{2i}\\
        &=\frac{e^{az+b_{1}+(b_{2}-ac-b_{1})}-e^{-az-b_{1}-(b_{2}-ac-b_{1})}}{2i}\\
        &=\frac{-aie^{az+b_{1}}-\frac{1}{-ai}e^{-az-b_{1}}}{2i}\\
        &=\frac{-ae^{az+b_{1}}-\frac{1}{a}e^{-az-b_{1}}}{2}.
\end{align*}
So we have
\begin{align*}
w'_{1}(z)=\frac{-a^2 e^{az+b_{1}}+e^{-az-b_{1}}}{2}.
\end{align*}
From \eqref{E:6.3} we also have
\begin{align*}
w'_{1}(z)=\frac{e^{az+b_{1}}+e^{-az-b_{1}}}{2},
\end{align*}
hence $a^2=-1$. That is $a=i$ or $a=-i$.

\par (1) If $a=i$, from \eqref{E:6.12a}, \eqref{E:6.12b} and $w_{1}(z+c)$, $w_{1}(z)$ as above, then
\begin{align}
e^{b_{2}-ic-b_{1}}=1, \label{E:6.12c}
\end{align}
\begin{align}
e^{b_{1}-ic-b_{2}}=1, \label{E:6.12d}
\end{align}
and
\begin{align*}
w_{1}(z+c)&=\frac{e^{iz+b_{2}}-e^{-iz-b_{2}}}{2i}\\
          &=\frac{e^{i(z-b_{2}i)}-e^{-i(z-b_{2}i)}}{2i}=\sin(z-b_{2}i),\\
w_{1}(z)&=\frac{-ie^{iz+b_{1}}-\frac{1}{i}e^{-iz-b_{1}}}{2}\\
        &=\frac{e^{iz+b_{1}}-e^{-iz-b_{1}}}{2i}=\sin(z-b_{1}i).
\end{align*}
From \eqref{E:6.4} and \eqref{E:6.12d} we also get
\begin{align*}
w_{2}(z+c)&=\frac{e^{iz+b_{1}}-e^{-iz-b_{1}}}{2i}\\
          &=\sin(z-b_{1}i),\\
w_{2}(z)&=\frac{e^{i(z-c)+b_{1}}-e^{-i(z-c)-b_{1}}}{2i}\\
        &=\frac{e^{iz+b_{2}+(b_{1}-ic-b_{2})}-e^{-iz-b_{2}-(b_{1}-ic-b_{2})}}{2i}\\
        &=\frac{e^{iz+b_{2}}-e^{-iz-b_{2}}}{2i}\\
        &=\sin(z-b_{2}i).
\end{align*}
Clearly $(w_{1}, w_{2})=(\sin(z-b_{1}i), \ \sin(z-b_{2}i))$ is a solution of Eq.\eqref{E:3.4}. From \eqref{E:6.12c} and \eqref{E:6.12d}, we have $e^{-2ic}=1$, $e^{2(b_{1}-b_{2})}=1$. That is $c=k\pi$ and $b_{1}-b_{2}=l\pi i$. Set
$B_{1}=-b_{1}i$, $B_{2}=-b_{2}i$, then $(w_{1}, w_{2})=(\sin(z+B_{1}),\ \sin(z+B_{2}))$, $B_{1}-B_{2}=l\pi$, $c=k\pi$,
$k,l \in \mathds{Z}$.

\par (2) If $a=-i$, similar as above, we have
\begin{align}
e^{b_{2}+ic-b_{1}}=-1, \label{E:6.12e}
\end{align}
\begin{align}
e^{b_{1}+ic-b_{2}}=-1, \label{E:6.12f}
\end{align}
and
\begin{align*}
w_{1}(z+c)&=\frac{e^{-iz+b_{2}}-e^{iz-b_{2}}}{2i}\\
          &=-\sin(z+b_{2}i),\\
w_{1}(z)&=\frac{ie^{-iz+b_{1}}+\frac{1}{i}e^{iz-b_{1}}}{2}\\
        &=\sin(z+b_{1}i).
\end{align*}
From \eqref{E:6.4} and \eqref{E:6.12f} we also get
\begin{align*}
w_{2}(z+c)&=\frac{e^{-iz+b_{1}}-e^{iz-b_{1}}}{2i}\\
          &=-\sin(z+b_{1}i),\\
w_{2}(z)&=\frac{e^{-i(z-c)+b_{1}}-e^{i(z-c)-b_{1}}}{2i}\\
        &=\frac{e^{-iz+b_{2}+(b_{1}+ic-b_{2})}-e^{iz-b_{2}-(b_{1}+ic-b_{2})}}{2i}\\
        &=\frac{e^{iz-b_{2}}-e^{-iz+b_{2}}}{2i}\\
        &=\sin(z+b_{2}i).
\end{align*}
We can easily verify that $(w_{1}, w_{2})=(\sin(z+b_{1}i), \ \sin(z+b_{2}i))$ is a solution of Eq.\eqref{E:3.4}. From \eqref{E:6.12e} and \eqref{E:6.12f}, we have $e^{2ic}=1$, $e^{2(b_{1}-b_{2})}=1$. That is $c=k\pi$ and $b_{1}-b_{2}=l\pi i$. Set
$B_{1}=b_{1}i$, $B_{2}=b_{2}i$, then $(w_{1}, w_{2})=(\sin(z+B_{1}),\ \sin(z+B_{2}))$, $B_{1}-B_{2}=l\pi$ ,$c=k\pi$,
$k,l \in \mathds{Z}$.

\par Subcase (ii): If $q=-1$, then $a_{1}=-a_{2}$. We can assume that $p(z)=az+b_{1}$ and $h(z)=-az+b_{2}$. Thus, by \eqref{E:6.11} and
\eqref{E:6.12}, we get
\begin{align}
\label{E:6.13}e^{b_{2}-ac-b_{1}}&=-ai,\\
\label{E:6.14}e^{b_{1}+ac-b_{2}}&=ai.
\end{align}
From \eqref{E:6.13} and \eqref{E:6.14}, we gain $a^{2}=1$. That is  $a=1$  or  $a=-1$.

\par (1) If $a=1$, from \eqref{E:6.13} or \eqref{E:6.14},  we get
\begin{align}
e^{b_{2}-c-b_{1}}&=-i,\label{E:6.14a}
\end{align}
and $p(z)=z+b_{1}$, $h(z)=-z+b_{2}$. From \eqref{E:6.6} and \eqref{E:6.14a}, we gain

\begin{align*}
w_{1}(-z+c)&=\frac{e^{-z+b_{2}}-e^{z-b_{2}}}{2i}\\
           &=i\sinh(z-b_{2}),\\
w_{1}(z)&=\frac{e^{z-c+b_{2}}-e^{-z+c-b_{2}}}{2i}\\
        &=\frac{e^{z+b_{1}+(b_{2}-c-b_{1})}-e^{-z-b_{1}-(b_{2}-c-b_{1})}}{2i}\\
        &=\frac{-ie^{z+b_{1}}-ie^{-z-b_{1}}}{2i}\\
        &=-\cosh{(z+b_{1})}.
\end{align*}
Similarly, from \eqref{E:6.4} and \eqref{E:6.14a}, we get
\begin{align*}
w_{2}(-z+c)&=\frac{e^{z+b_{1}}-e^{-z-b_{1}}}{2i}\\
           &=-i\sinh(z+b_{1}),\\
w_{2}(z)&=\frac{e^{-z+c+b_{1}}-e^{z-c-b_{1}}}{2i}\\
&=\frac{e^{-z+b_{2}+(b_{1}+c-b_{2})}-e^{z-b_{2}-(b_{1}+c-b_{2})}}{2i}\\
&=\frac{ie^{-z+b_{2}}-(-i)e^{z-b_{2}}}{2i}\\
&=\cosh(z-b_{2}).
\end{align*}

By the property of hyperbolic functions, we get $w^{\prime}_{1}(z)=-\sinh{(z+b_{1})}$, $w'_{2}(z)=\sinh(z-b_{2})$. It is easy to see that $(w_{1}(z),w_{2}(z))$ is not a solution of Eq.\eqref{E:3.4}.

\par (2) If $a=-1$, from \eqref{E:6.13} or \eqref{E:6.14},we have
\begin{align}
e^{b_{2}+c-b_{1}}&=i,\label{E:6.14b}
\end{align}
and $p(z)=-z+b_{1}$, $h(z)=z+b_{2}$. From \eqref{E:6.6}, \eqref{E:6.4} and \eqref{E:6.14b} we gain
\begin{align*}
w_{1}(-z+c)&=\frac{e^{z+b_{2}}-e^{-z-b_{2}}}{2i}\\
           &=-i\sinh(z+b_{2}),\\
w_{1}(z)&=\frac{e^{-z+c+b_{2}}-e^{z-c-b_{2}}}{2i}\\
        &=\frac{e^{-z+b_{1}+(b_{2}+c-b_{1})}-e^{z-b_{1}-(b_{2}+c-b_{1})}}{2i}\\
        &=\frac{ie^{-z+b_{1}}+ie^{z-b_{1}}}{2i}\\
        &=\cosh(z-b_{1}),
\end{align*}
and
\begin{align*}
w_{2}(-z+c)&=\frac{e^{-z+b_{1}}-e^{z-b_{1}}}{2i}\\
           &=i\sinh(z-b_{1}),\\
w_{2}(z)&=\frac{e^{z-c+b_{1}}-e^{-z+c-b_{1}}}{2i}\\
&=\frac{e^{z+b_{2}-(b_{2}+c-b_{1})}-e^{-z-b_{2}+(b_{2}+c-b_{1})}}{2i}\\
&=\frac{-ie^{z+b_{2}}-ie^{-z-b_{2}}}{2i}\\
&=-\cosh(z+b_{2}).
\end{align*}
Obviously $(w_{1}(z),w_{2}(z))$ is not a solution of Eq.\eqref{E:3.4}.

Hence, in this case, if $q=-1$, Eq.\eqref{E:3.4} does not posses any finite order of transcendental entire solutions.

\par {\bf Case (B):} Let $e^{2p(z)}$, $e^{2h(z)}$, $\frac{q i}{h^{\prime}(z)}e^{h(z)-p(qz+c)}$ and $\frac{q i}{p^{\prime}(z)}e^{p(z)-h(qz+c)}$ be not constants. From \eqref{E:6.9}, \eqref{E:6.10} and Lemma 4.2, we get
\begin{align}\label{E:6.18}
\frac{q i}{h^{\prime}(z)}e^{h(z)+p(qz+c)}\equiv 1,
\end{align}
and
\begin{align}\label{E:6.19}
\frac{q i}{p^{\prime}(z)}e^{p(z)+h(qz+c)}\equiv 1.
\end{align}
Similar as above, $p(z)$ and $h(z)$ are nonconstant polynomials of first order.
Let $p(z)=a_{1}z+b_{1}$ and $h(z)=a_{2}z+b_{2}$, where $a_{1}(\neq0)$, $a_{2}(\neq0)$, $b_{1}$ and $b_{2}$ are constants. From \eqref{E:6.18}, we get
$$\frac{q i}{h^{\prime}(z)}e^{h(z)+p(qz+c)}=\frac{q i}{a_{2}}e^{a_{2}z+b_{2}+a_{1}(qz+c)+b_{1}}\equiv 1.$$
Thus, we can deduce that $a_{2}+qa_{1}=0$. Also from \eqref{E:6.19}, we can gain $a_{1}+qa_{2}=0$.  Then $a_{1}a_{2}=q^{2}a_{1}a_{2}$. Moreover, $q=1$ or $q=-1$.

\par Subcase (i): If $q=1$, then $a_{2}= -a_{1}$. Thus, we assume that $p(z)=az+b_{1}$ and $h(z)=-az+b_{2}$, where $a(\neq0)$, $b_{1}$ and $b_{2}$ are constants.
From \eqref{E:6.18} and \eqref{E:6.19}, we obtain
\begin{align}
\label{E:6.20}e^{b_{1}+b_{2}+ac}=ai
\end{align}
and
\begin{align}
\label{E:6.21}e^{b_{1}+b_{2}-ac}=-ai.
\end{align}
From \eqref{E:6.3}, we have
\begin{align}
\label{E:6.22}w'_{1}(z)=\frac{e^{az+b_{1}}+e^{-az-b_{1}}}{2}.
\end{align}
From \eqref{E:6.6} and \eqref{E:6.20}, we gain
\begin{align}
\label{E:6.22a}w_{1}(z+c)&=\frac{e^{-az+b_{2}}-e^{az-b_{2}}}{2i},\\
\nonumber w_{1}(z)&=\frac{e^{-a(z-c)+b_{2}}-e^{a(z-c)-b_{2}}}{2i}\\
\nonumber&=\frac{e^{-az-b_{1}+(b_{1}+b_{2}+ac)}-e^{az+b_{1}-(b_{1}+b_{2}+ac)}}{2i}\\
\nonumber&=\frac{aie^{-az-b_{1}}-\frac{1}{ai}e^{az+b_{1}}}{2i}\\
\label{E:6.23}&=\frac{ae^{-az-b_{1}}+\frac{1}{a}e^{az+b_{1}}}{2},
\end{align}
and
\begin{align*}
 w'_{1}(z)=\frac{-a^2 e^{-az-b_{1}}+e^{az+b_{1}}}{2}.
\end{align*}
Combining with \eqref{E:6.22}, then $a^2=-1$. That is $a=i$ or $a=-i$.

\par (1) If $a=i$, by \eqref{E:6.20}, \eqref{E:6.21}, we get
\begin{align}
e^{b_{1}+b_{2}+ic}=-1, \label{E:6.21a}\\
e^{b_{1}+b_{2}-ic}=1, \label{E:6.21b}
\end{align}
From \eqref{E:6.22a}, \eqref{E:6.23}, we gain
\begin{align*}
w_{1}(z+c)&=\frac{e^{-iz+b_{2}}-e^{iz-b_{2}}}{2i}\\
          &=-\sin{(z+b_{2}i)},\\
w_{1}(z)&=\frac{ie^{-iz-b_{1}}+\frac{1}{i}e^{iz+b_{1}}}{2}\\
        &=\frac{e^{iz+b_{1}}-e^{-iz-b_{1}}}{2i}\\
        &=\sin{(z-b_{1}i)}.
\end{align*}
From \eqref{E:6.4} and \eqref{E:6.21b}, we get
\begin{align*}
w_{2}(z+c)&=\frac{e^{iz+b_{1}}-e^{-iz-b_{1}}}{2i}\\
          &=\sin(z-b_{1}i),\\
w_{2}(z)&=\frac{e^{i(z-c)+b_{1}}-e^{-i(z-c)-b_{1}}}{2i}\\
        &=\frac{e^{iz-b_{2}+(b_{1}+b_{2}-ci)}-e^{-iz+b_{2}-(b_{1}+b_{2}-ci)}}{2i}\\
        &=\frac{e^{iz-b_{2}}-e^{-iz+b_{2}}}{2i}\\
        &=\sin(z+b_{2}i).
\end{align*}
It is easy to be verified that  $(w_{1}, w_{2})=(\sin(z-b_{1}i),\ \sin(z+b_{2}i))$ is a solution of Eq.\eqref{E:3.4}.
From \eqref{E:6.21a}, \eqref{E:6.21b}, we deduce that $e^{2ci}=-1$, $2(b_{1}+b_{2})=-1$. So $c=\frac{\pi}{2}+k\pi$,
$b_{1}+b_{2}=\frac{\pi}{2}i+k\pi i$. Set $B_{1}=-b_{1}i$, $B_{2}=b_{2}i$, then
$(w_{1}, w_{2})=(\sin(z+B_{1}),\ \sin(z+B_{2}))$, $B_{1}-B_{2}=\frac{\pi}{2}+l\pi$, $c=\frac{\pi}{2}+k\pi$,
$k, l \in \mathds{Z}$. This is the result we want.

\par(2) If $a=-i$, from \eqref{E:6.20}, \eqref{E:6.21}, we know \eqref{E:6.21a} and \eqref{E:6.21b} still hold.
From \eqref{E:6.22a}, \eqref{E:6.23}, we have
\begin{align*}
w_{1}(z+c)&=\frac{e^{iz+b_{2}}-e^{-iz-b_{2}}}{2i}\\
          &=\sin{(z-b_{2}i)},\\
w_{1}(z)&=\frac{-ie^{iz-b_{1}}+\frac{1}{-i}e^{-iz+b_{1}}}{2}\\
        &=\frac{e^{iz-b_{1}}-e^{-iz+b_{1}}}{2i}\\
        &=\sin{(z+b_{1}i)}.
\end{align*}
From \eqref{E:6.4}and \eqref{E:6.21b}, we have
\begin{align*}
w_{2}(z+c)&=\frac{e^{-iz+b_{1}}-e^{iz-b_{1}}}{2i}\\
          &=-\sin(z+b_{1}i),\\
w_{2}(z)&=\frac{e^{-i(z-c)+b_{1}}-e^{i(z-c)-b_{1}}}{2i}\\
        &=\frac{e^{-iz-b_{2}+(b_{1}+b_{2}+ci)}-e^{iz+b_{2}-(b_{1}+b_{2}+ci)}}{2i}\\
        &=\frac{e^{iz+b_{2}}-e^{-iz-b_{2}}}{2i}\\
        &=\sin(z-b_{2}i).
\end{align*}
It is easy to see that $(w_{1}, w_{2})=(\sin(z+b_{1}i),\ \sin(z-b_{2}i))$ is a solution of Eq.\eqref{E:3.4}.
Set $B_{1}=-b_{1}i$, $B_{2}=b_{2}i$, from \eqref{E:6.21a}, \eqref{E:6.21b}, then
$(w_{1}, w_{2})=(\sin(z+B_{1}),\ \sin(z+B_{2}))$,
$B_{1}-B_{2}=\frac{\pi}{2}+l\pi$, $c=\frac{\pi}{2}+k\pi$, $k, l \in \mathds{Z}$.

\par Subcase (ii): If $q=-1$, then $a_{2}=a_{1}$, so we can assume that $p(z)=az+b_{1}$ and $h(z)=az+b_{2}$.
Thus, by \eqref{E:6.18} or \eqref{E:6.19}, we get
Then
\begin{align}\label{E:6.30}e^{b_{2}+ac+b_{1}}=ai.
\end{align}
From \eqref{E:6.3}, we get
\begin{align}
\label{E:6.32}w'_{1}(z)=\frac{e^{az+b_{1}}+e^{-az-b_{1}}}{2},
\end{align}
By \eqref{E:6.6} and \eqref{E:6.30}, we have
 \begin{align}
\label{E:6.33} w_{1}(-z+c)&=\frac{e^{az+b_{2}}-e^{-az-b_{2}}}{2i},\\
\nonumber         w_{1}(z)&=\frac{e^{a(-z+c)+b_{2}}-e^{-a(-z+c)-b_{2}}}{2i}\\
\nonumber                 &=\frac{e^{-az-b_{1}+(b_{1}+b_{2}+ac)}-e^{az+b_{1}-(b_{1}+b_{2}+ac)}}{2i}\\
\nonumber                 &=\frac{aie^{-az-b_{1}}-\frac{1}{ai}e^{az+b_{1}}}{2i}\\
\label{E:6.34}            &=\frac{ae^{-az-b_{1}}+\frac{1}{a}e^{az+b_{1}}}{2},
\end{align}
and
\begin{align}
\nonumber        w'_{1}(z)&=\frac{-a^2e^{-az-b_{1}}+e^{az+b_{1}}}{2}.
\end{align}
Comparing with \eqref{E:6.32}, we gain $a^2=-1$. That is  $a=i$  or  $a=-i$.

\par (1) If $a=i$, from \eqref{E:6.30}, we get
\begin{align}
\label{E:6.30a} e^{b_{2}+ic+b_{1}}=-1.
\end{align}
By \eqref{E:6.33}, \eqref{E:6.34} we obtain
\begin{align*}
w_{1}(-z+c)&=\frac{e^{iz+b_{2}}-e^{-iz-b_{2}}}{2i}\\
           &=\sin{(z-b_{2}i)},\\
w_{1}(z)&=\frac{ie^{-iz-b_{1}}+\frac{1}{i}e^{iz+b_{1}}}{2}\\
        &=\sin(z-b_{1}i)
\end{align*}
From \eqref{E:6.4} and \eqref{E:6.30a}, we get
\begin{align*}
w_{2}(-z+c)&=\frac{e^{iz+b_{1}}-e^{-iz-b_{1}}}{2i}\\
           &=\sin{(z-b_{1}i)},\\
w_{2}(z)&=\frac{e^{i(-z+c)+b_{1}}-e^{-i(-z+c)-b_{1}}}{2i}\\
&=\frac{e^{-iz-b_{2}+(b_{1}+b_{2}+ic)}-e^{iz+b_{2}-(b_{1}+b_{2}+ic)}}{2i}\\
&=\frac{-e^{-iz-b_{2}}+e^{iz+b_{2}}}{2i}\\
&=\sin{(z-b_{2}i)}.
\end{align*}
Obviously $(w_{1}, w_{2})=(\sin(z-b_{1}i),\ \sin(z-b_{2}i))$ satisfies Eq.\eqref{E:3.4}. By \eqref{E:6.30a}, we know $b_{1}+b_{2}+ic=(2k+1)\pi i~~(k\in \mathds{Z})$. Set $B_{1}=-b_{1}i$ and $B_{2}=-b_{2}i$, then
$(w_{1}, w_{2})=(\sin(z+B_{1}),\ \sin(z+B_{2}))$, $B_{1}+B_{2}+c=(2k+1)\pi$, $k\in \mathds{Z}$.

\par (2) If $a=-i$, from \eqref{E:6.30}, we derive that
\begin{align}
\label{E:6.30b} e^{b_{1}-ic+b_{2}}=1.
\end{align}
From \eqref{E:6.33}, \eqref{E:6.34}, we get
\begin{align*}
w_{1}(-z+c)&=\frac{e^{-iz+b_{2}}-e^{iz-b_{2}}}{2i}\\
           &=-\sin{(z+b_{2}i)},\\
w_{1}(z)&=\frac{-ie^{iz-b_{1}}+\frac{1}{-i}e^{-iz+b_{1}}}{2}\\
        &=\sin(z+b_{1}i).
\end{align*}
From \eqref{E:6.4} and \eqref{E:6.30b}, we gain
\begin{align*}
w_{2}(-z+c)&=\frac{e^{-iz+b_{1}}-e^{iz-b_{1}}}{2i}\\
           &=-\sin{(z+b_{1}i)},\\
w_{2}(z)&=\frac{e^{-i(-z+c)+b_{1}}-e^{i(-z+c)-b_{1}}}{2i}\\
&=\frac{e^{iz-b_{2}+(b_{1}+b_{2}-ic)}-e^{-iz+b_{2}-(b_{1}+b_{2}-ic)}}{2i}\\
&=\frac{e^{iz-b_{2}}-e^{-iz+b_{2}}}{2i}\\
&=\sin{(z+b_{2}i)}.
\end{align*}
It is clear that $(w_{1}, w_{2})=(\sin(z+b_{1}i),\ \sin(z+b_{2}i))$ is a solution of Eq.\eqref{E:3.4}.
From \eqref{E:6.30b} we get $b_{2}+b_{1}-ic=2k\pi i~~(k\in \mathds{Z})$. Set $B_{1}=b_{1}i$ and $B_{2}=b_{2}i$,
then $(w_{1}, w_{2})=(\sin(z+B_{1}),\ \sin(z+B_{2}))$, $B_{1}+B_{2}+c=2k\pi$, $k\in \mathds{Z}$.

\par {\bf Case (C):} Let $e^{2p(z)}$, $e^{2h(z)}$, $\frac{q i}{h^{\prime}(z)}e^{h(z)+p(qz+c)}$ and $\frac{q i}{p^{\prime}(z)}e^{p(z)-h(qz+c)}$ be not constants, From \eqref{E:6.9}, \eqref{E:6.10} and Lemma 4.2, we get
\begin{align}\label{E:6.38}
\frac{q i}{h^{\prime}(z)}e^{h(z)-p(qz+c)}\equiv 1,
\end{align}
and
\begin{align}\label{E:6.39}
\frac{q i}{p^{\prime}(z)}e^{p(z)+h(qz+c)}\equiv 1.
\end{align}
As same as above, $p(z)$ and $h(z)$ are first order polynomials.
Let $p(z)=a_{1}z+b_{1}$ and $h(z)=a_{2}z+b_{2}$, where $a_{1}(\neq0)$, $a_{2}(\neq0)$, $b_{1}$ and $b_{2}$ are constants. From \eqref{E:6.38} and \eqref{E:6.39}, we get
$a_{2}-qa_{1}=0$ and $a_{1}+qa_{2}=0$.  Then $a_{1}a_{2}=-q^{2}a_{1}a_{2}$. Furthermore, we get $q=i$ or $q=-i$.

\par Subcase (i): If $q=i$, then $a_{2}= ia_{1}$. Thus, we assume that $p(z)=az+b_{1}$ and $h(z)=iaz+b_{2}$, where $a(\neq0)$, $b_{1}$ and $b_{2}$ are constants.
From \eqref{E:6.38} and \eqref{E:6.39}, we obtain
\begin{align}
\label{E:6.40}e^{b_{2}-b_{1}-ac}=-ai,
\end{align}
and
\begin{align}
\label{E:6.41}e^{b_{1}+b_{2}+iac}=-a.
\end{align}
From \eqref{E:6.3}, we get
\begin{align}
\label{E:6.42}w'_{1}(z)=\frac{e^{az+b_{1}}+e^{-az-b_{1}}}{2}.
\end{align}
By \eqref{E:6.6} and \eqref{E:6.40}, we gain
 \begin{align}
\label{E:6.42a} w_{1}(iz+c)&=\frac{e^{iaz+b_{2}}-e^{-iaz-b_{2}}}{2i},\\
\nonumber w_{1}(z)&=\frac{e^{ia(-iz+ic)+b_{2}}-e^{-ia(-iz+ic)-b_{2}}}{2i}\\
\nonumber&=\frac{e^{az+b_{1}+(b_{2}-b_{1}-ac)}-e^{-az-b_{1}-(b_{2}-b_{1}-ac)}}{2i}\\
\nonumber&=\frac{-aie^{az+b_{1}}-\frac{1}{-ai}e^{-az-b_{1}}}{2i}\\
\label{E:6.43}&=\frac{-ae^{az+b_{1}}-\frac{1}{a}e^{-az-b_{1}}}{2},
\end{align}
and
\begin{align}
\nonumber w'_{1}(z)=\frac{-a^2e^{az+b_{1}}+e^{-az-b_{1}}}{2}.
\end{align}
Comparing with \eqref{E:6.42}, we get $-a^{2}=1$, then $a=i$ or $a=-i$.

\par (1) If $a=i$, then $p(z)=iz+b_{1}$ and $h(z)=-z+b_{2}$. By \eqref{E:6.40}, \eqref{E:6.41}, we get
\begin{align}
\label{E:6.41a} e^{b_{2}-b_{1}-ic}&=1,\\
\label{E:6.41b} e^{b_{1}+b_{2}-c}&=-i.
\end{align}
From \eqref{E:6.42a}, \eqref{E:6.43} we gain
\begin{align*}
w_{1}(iz+c)&=\frac{e^{-z+b_{2}}-e^{z-b_{2}}}{2i}\\
           &=i\sinh(z-b_{2}),\\
w_{1}(z)&=\frac{-ie^{iz+b_{1}}-\frac{1}{i}e^{-iz-b_{1}}}{2}\\
        &=\sin(z-b_{1}i).
\end{align*}
From \eqref{E:6.4} and \eqref{E:6.41b}, we gain
\begin{align*}
\nonumber w_{2}(iz+c)&=\frac{e^{iz+b_{1}}-e^{-iz-b_{1}}}{2i}\\
                     &=\sin(z-b_{1}i),\\
 w_{2}(z)&=\frac{e^{i(-iz+ic)+b_{1}}-e^{-i(-zi+ci)-b_{1}}}{2i}\\
 &=\frac{e^{z-b_{2}+(b_{1}+b_{2}-c)}-e^{-z+b_{2}-(b_{1}+b_{2}-c)}}{2i}\\
 &=\frac{-ie^{z-b_{2}}-\frac{1}{-i}e^{-z+b_{2}}}{2i}\\
 &=-\cosh(z-b_{2}).
\end{align*}
Obviously $(w_{1},w_{2})=(\sin(z-b_{1}i),-\cosh(z-b_{2}))$ is not a solution of Eq.\eqref{E:3.4}.

\par (2) If $a=-i$, we assume that $p(z)=-iz+b_{1}$ and $h(z)=z+b_{2}$. By \eqref{E:6.41}, \eqref{E:6.41},
we get
\begin{align}
\label{E:6.41c} e^{b_{2}-b_{1}+ic}=1,\\
\label{E:6.41d} e^{b_{1}+b_{2}+c}=i.
\end{align}
From \eqref{E:6.42a}, \eqref{E:6.43}, we have
\begin{align*}
 w_{1}(iz+c)&=\frac{e^{z+b_{2}}-e^{-z-b_{2}}}{2i}\\
            &=-i\sinh(z+b_{2}),\\
 w_{1}(z)&=\frac{ie^{-iz+b_{1}}-\frac{1}{-i}e^{iz-b_{1}}}{2}\\
         &=\sin{(z+b_{1}i)}.
\end{align*}
From \eqref{E:6.4} and \eqref{E:6.41d}, we obtain
\begin{align*}
w_{2}(iz+c)&=\frac{e^{-iz+b_{1}}-e^{iz-b_{1}}}{2i}\\
           &=-\sin(z+b_{1}i),\\
w_{2}(z)&=\frac{e^{-i(-iz+ic)+b_{1}}-e^{i(-zi+ci)-b_{1}}}{2i}\\
&=\frac{e^{-z-b_{2}+(b_{1}+b_{2}+c)}-e^{z+b_{2}-(b_{1}+b_{2}+c)}}{2i}\\
&=\frac{ie^{-z-b_{2}}-\frac{1}{i}e^{z+b_{2}}}{2i}\\
&=\cosh(z+b_{2}).
\end{align*}
It is easy to see that $(w_{1},w_{2})=(\sin(z+b_{1}i),\cosh(z+b_{2}))$ is not a solution of Eq.\eqref{E:3.4}.

\par Subcase (ii): If $q=-i$, then $a_{2}=-ia_{1}$. Thus, let $p(z)=az+b_{1}$ and $h(z)=-iaz+b_{2}$, where $a(\neq 0)$, $b_{1}$ and $b_{2}$ are constants. From \eqref{E:6.38} and \eqref{E:6.39}, we obtain
\begin{align}
\label{E:6.49}e^{b_{2}-b_{1}-ac}&=-ai,\\
\label{E:6.50}e^{b_{1}+b_{2}-iac}&=a.
\end{align}
From \eqref{E:6.3}, we get
\begin{align}
\label{E:6.51}w'_{1}(z)=\frac{e^{az+b_{1}}+e^{-az-b_{1}}}{2}.
\end{align}
By \eqref{E:6.6}, \eqref{E:6.49} we have
\begin{align}
\label{E:6.51a} w_{1}(-iz+c)&=\frac{e^{-iaz+b_{2}}-e^{iaz-b_{2}}}{2i},\\
\nonumber w_{1}(z)&=\frac{e^{-ia(iz-ic)+b_{2}}-e^{ia(iz-ic)-b_{2}}}{2i}\\
\nonumber&=\frac{e^{az+b_{1}+(b_{2}-b_{1}-ac)}-e^{-az-b_{1}-(b_{2}-b_{1}-ac)}}{2i}\\
\nonumber&=\frac{-aie^{az+b_{1}}-\frac{1}{-ai}e^{-az-b_{1}}}{2i}\\
\label{E:6.52}&=\frac{-ae^{az+b_{1}}-\frac{1}{a}e^{-az-b_{1}}}{2},
\end{align}
and
\begin{align}
\nonumber w'_{1}(z)=\frac{-a^{2}e^{az+b_{1}}+e^{-az-b_{1}}}{2}.
\end{align}
Comparing with \eqref{E:6.51}, we gain $-a^{2}=-1$. That is  $a=i$  or  $a=-i$.

\par (1) If $a=i$, we can assume $p(z)=iz+b_{1}$, $h(z)=z+b_{2}$. From \eqref{E:6.49}, \eqref{E:6.50}, we get
\begin{align}
\label{E:6.50a} e^{b_{2}-b_{1}-ic}&=1.\\
\label{E:6.50b} e^{b_{2}+b_{1}+c}&=i.
\end{align}
By \eqref{E:6.51a}, \eqref{E:6.52} we obtain
\begin{align*}
w_{1}(-iz+c)&=\frac{e^{z+b_{2}}-e^{-iz-b_{2}}}{2i}\\
            &=-i\sinh(z+b_{2}),\\
w_{1}(z)&=\frac{-ie^{iz+b_{1}}-\frac{1}{i}e^{-iz-b_{1}}}{2}\\
        &=\sin(z-b_{1}i).
\end{align*}
From \eqref{E:6.4}, \eqref{E:6.50b} we also get
\begin{align*}
w_{2}(-iz+c)&=\frac{e^{iz+b_{1}}-e^{-iz-b_{1}}}{2i}\\
            &=\sin(z-b_{1}i),\\
w_{2}(z)&=\frac{e^{i(iz-ic)+b_{1}}-e^{-i(iz-ic)-b_{1}}}{2i}\\
&=\frac{e^{-z-b_{2}+(b_{2}+b_{1}+c)}-e^{z+b_{2}-(b_{2}+b_{1}+c)}}{2i}\\
&=\frac{ie^{-z-b_{2}}-\frac{1}{i}e^{z+b_{2}}}{2i}\\
&=\cosh{(z+b_{2})}.
\end{align*}
It is clear that $(w_{1}, w_{2})=(\sin(z-b_{1}i), \cosh(z+b_{2}))$ is not a solution of Eq.\eqref{E:3.4}

\par (2) If $a=-i$, we can set $p(z)=-iz+b_{1}$ and $h(z)=-z+b_{2}$. From \eqref{E:6.49}, \eqref{E:6.50}, we get
\begin{align}
\label{E:6.50c} e^{b_{2}-b_{1}+ic}&=-1,
\end{align}
and
\begin{align}
\label{E:6.50d} e^{b_{1}+b_{2}-c}&=-i.
\end{align}
Thus, by \eqref{E:6.51a}, \eqref{E:6.52}, we get
\begin{align*}
 w_{1}(-iz+c)&=\frac{e^{-z+b_{2}}-e^{z-b_{2}}}{2i}\\
             &=i\sinh(z-b_{2}),\\
 w_{1}(z)&=\frac{ie^{-iz+b_{1}}-\frac{1}{-i}e^{iz-b_{1}}}{2}\\
         &=\sin{(z+b_{1}i)}.
\end{align*}
From \eqref{E:6.4} and  \eqref{E:6.50d}, we deduce that
\begin{align*}
 w_{2}(-iz+c)&=\frac{e^{-iz+b_{1}}-e^{iz-b_{1}}}{2i}\\
             &=-i\sin(z+b_{1}i),\\
 w_{2}(z)&=\frac{e^{-i(iz-ic)+b_{1}}-e^{i(iz-ic)-b_{1}}}{2i}\\
 &=\frac{e^{z-b_{2}+(b_{1}+b_{2}-c)}-e^{-z+b_{2}-(b_{1}+b_{2}-c)}}{2i}\\
 &=\frac{-ie^{z-b_{2}}-\frac{1}{-i}e^{-z+b_{2}}}{2i}\\
 &=-\cosh(z-b_{2}).
\end{align*}
It is easy verification that $(w_{1}, w_{2})=(\sin(z+b_{1}i),\ -\cosh(z-b_{2}))$ is not a solution of Eq.\eqref{E:3.4}.

\par {\bf Case (D):} Let $e^{2p(z)}$, $e^{2h(z)}$, $\frac{q i}{h^{\prime}(z)}e^{h(z)-p(qz+c)}$ and $\frac{q i}{p^{\prime}(z)}e^{p(z)+h(qz+c)}$ be not constants, From \eqref{E:6.9}, \eqref{E:6.10} and Lemma 4.2, we get
\begin{align}\label{E:6.56}
\frac{q i}{h^{\prime}(z)}e^{h(z)+p(qz+c)}\equiv 1,
\end{align}
and
\begin{align}\label{E:6.57}
\frac{q i}{p^{\prime}(z)}e^{p(z)-h(qz+c)}\equiv 1.
\end{align}
As similar as above $p(z)$ and $h(z)$ are first order polynomials.
Hence, Let $p(z)=a_{1}z+b_{1}$ and $h(z)=a_{2}z+b_{2}$, where $a_{1}(\neq0)$, $a_{2}(\neq0)$, $b_{1}$ and $b_{2}$ are constants. From \eqref{E:6.56}, \eqref{E:6.57} we get
$a_{2}+qa_{1}=0$ and $a_{1}-qa_{2}=0$.  Then $a_{1}a_{2}=-q^{2}a_{1}a_{2}$. That is $q=i$ or $q=-i$.

\par Subcase (i): If $q=i$, then $a_{2}=-ia_{1}$. Therefore, let $p(z)=az+b_{1}$ and $h(z)=-iaz+b_{2}$, where $a(\neq0)$, $b_{1}$ and $b_{2}$ are constants.
From \eqref{E:6.56} and \eqref{E:6.57}, we obtain
\begin{align}
\label{E:6.58}e^{b_{1}+b_{2}+ac}&=ai\\
\label{E:6.59}e^{b_{1}-b_{2}+iac}&=-a.
\end{align}
From \eqref{E:6.3}, we get
\begin{align}
\label{E:6.60}w'_{1}(z)=\frac{e^{az+b_{1}}+e^{-az-b_{1}}}{2}.
\end{align}
By \eqref{E:6.6} and \eqref{E:6.58}, we gain
\begin{align}
\label{E:6.60a}w_{1}(iz+c)&=\frac{e^{-iaz+b_{2}}-e^{iaz-b_{2}}}{2i},\\
\nonumber w_{1}(z)&=\frac{e^{-ia(-iz+ic)+b_{2}}-e^{ia(-iz+ic)-b_{2}}}{2i}\\
\nonumber&=\frac{e^{-az-b_{1}+(b_{1}+b_{2}+ac)}-e^{az+b_{1}-(b_{1}+b_{2}+ac)}}{2i}\\
\nonumber&=\frac{aie^{-az-b_{1}}-\frac{1}{ai}e^{az+b_{1}}}{2i}\\
\label{E:6.61}&=\frac{ae^{-az-b_{1}}+\frac{1}{a}e^{az+b_{1}}}{2},
\end{align}
and
\begin{align}
\nonumber w'_{1}(z)=\frac{-a^{2}e^{-az-b_{1}}+e^{az+b_{1}}}{2}.
\end{align}
Comparing with \eqref{E:6.60}, we get $a^{2}=-1$. That is $a=i$ or $a=-i$.

\par (1) If $a=i$, we may assume that $p(z)=iz+b_{1}$, $h(z)=z+b_{2}$. By \eqref{E:6.58}, \eqref{E:6.59}, we get
\begin{align}
\label{E:6.59a} e^{b_{1}+b_{2}+ic}&=-1,\\
\label{E:6.59b} e^{b_{1}-b_{2}-c}&=-i.
\end{align}
From \eqref{E:6.60a}, \eqref{E:6.61}, we gain
\begin{align*}
w_{1}(iz+c)&=\frac{e^{z+b_{2}}-e^{-z-b_{2}}}{2i}\\
           &=-i\sinh(z+b_{2}),\\
w_{1}(z)&=\frac{ie^{-iz-b_{1}}+\frac{1}{i}e^{iz+b_{1}}}{2}\\
        &=\sin{(z-b_{1}i)}.
\end{align*}
From \eqref{E:6.4} and \eqref{E:6.59b}, we gain
\begin{align*}
w_{2}(iz+c)&=\frac{e^{iz+b_{1}}-e^{-iz-b_{1}}}{2i}\\
           &=\sin(z-b_{1}),\\
w_{2}(z)&=\frac{e^{i(-iz+ic)+b_{1}}-e^{-i(-zi+ci)-b_{1}}}{2i}\\
&=\frac{e^{z+b_{2}+(b_{1}-b_{2}-c)}-e^{-z-b_{2}-(b_{1}-b_{2}-c)}}{2i}\\
&=\frac{-ie^{z+b_{2}}-\frac{1}{-i}e^{-z-b_{2}}}{2i}\\
&=-\cosh(z+b_{2}).
\end{align*}
Obviously, $(w_{1},w_{2})=(\sin(z-b_{1}i),-\cosh(z+b_{2}))$ is not a solution of Eq.\eqref{E:3.4}.

\par (2) If $a=-i$, we can assume that $p(z)=-iz+b_{1}$ and $h(z)=-z+b_{2}$. By \eqref{E:6.58} and \eqref{E:6.59}, we get
\begin{align}
\label{E:6.59c} e^{b_{1}+b_{2}-ic}&=1,\\
\label{E:6.59d} e^{b_{1}-b_{2}+c}&=i.
\end{align}
From \eqref{E:6.60a}, \eqref{E:6.61}, we gain
\begin{align*}
 w_{1}(iz+c)&=\frac{e^{-z+b_{2}}-e^{z-b_{2}}}{2i}\\
            &=i\sinh(z-b_{2}),\\
 w_{1}(z)&=\frac{-ie^{iz-b_{1}}+\frac{1}{-i}e^{-iz+b_{1}}}{2}\\
         &=\sin{(z+b_{1}i)}.
\end{align*}
From \eqref{E:6.4} and \eqref{E:6.59d}, we gain
\begin{align*}
w_{2}(iz+c)&=\frac{e^{-iz+b_{1}}-e^{iz-b_{1}}}{2i}\\
           &=-\sin(z+b_{1}i),\\
w_{2}(z)&=\frac{e^{-i(-iz+ic)+b_{1}}-e^{i(-zi+ci)-b_{1}}}{2i}\\
        &=\frac{e^{-z+b_{2}+(b_{1}-b_{2}+c)}-e^{z-b_{2}-(b_{1}-b_{2}+c)}}{2i}\\
        &=\frac{ie^{-z+b_{2}}-\frac{1}{i}e^{z-b_{2}}}{2i}\\
        &=\cosh{(z-b_{2})}.
\end{align*}
It is easy to be verified that $(w_{1},w_{2})=(\sin(z+b_{1}i),\cosh(z-b_{2}))$ is not a solution of Eq.\eqref{E:3.4}.

\par Subcase (ii): If $q=-i$, then $a_{2}=ia_{1}$. Thus, we assume that $p(z)=az+b_{1}$ and $h(z)=iaz+b_{2}$, where $a(\neq 0)$, $b_{1}$ and $b_{2}$ are constants. Thus, from \eqref{E:6.56}, \eqref{E:6.56}, we get
\begin{align}
\label{E:6.66}e^{b_{1}+b_{2}+ac}&=ai,\\
\label{E:6.67}e^{b_{1}-b_{2}-iac}&=a.
\end{align}
From \eqref{E:6.3}, we get
\begin{align}
\label{E:6.68}w'_{1}(z)=\frac{e^{az+b_{1}}+e^{-az-b_{1}}}{2}.
\end{align}
By \eqref{E:6.6} and \eqref{E:6.66}, we have
\begin{align}
\label{E:6.68a} w_{1}(-iz+c)&=\frac{e^{iaz+b_{2}}-e^{-iaz-b_{2}}}{2i},\\
\nonumber w_{1}(z)&=\frac{e^{ia(iz-ic)+b_{2}}-e^{-ia(iz-ic)-b_{2}}}{2i}\\
\nonumber&=\frac{e^{-az-b_{1}+(b_{1}+b_{2}+ac)}-e^{az+b_{1}-(b_{1}+b_{2}+ac)}}{2i}\\
\nonumber&=\frac{aie^{-az-b_{1}}-\frac{1}{ai}e^{az+b_{1}}}{2i}\\
\label{E:6.69}&=\frac{ae^{-az-b_{1}}+\frac{1}{a}e^{az+b_{1}}}{2}.
\end{align}
and
\begin{align}
\nonumber w'_{1}(z)=\frac{-a^{2}e^{-az-b_{1}}+e^{az+b_{1}}}{2}.
\end{align}
Comparing with \eqref{E:6.68}, we gain $a^{2}=-1$. We conclude that  $a=i$  or  $a=-i$.

\par (1) If $a=i$, we can take $p(z)=iz+b_{1}$, $h(z)=-z+b_{2}$. From \eqref{E:6.66} and \eqref{E:6.67}, we get
\begin{align}
\label{E:6.67a} e^{b_{1}+b_{2}+ic}&=-1, \\
\label{E:6.67b} e^{b_{1}-b_{2}+c}&=i.
\end{align}
By \eqref{E:6.68a} and \eqref{E:6.69}, we obtain
\begin{align*}
w_{1}(-iz+c)&=\frac{e^{-z+b_{2}}-e^{z-b_{2}}}{2i}\\
            &=i\sinh(z-b_{2}),\\
w_{1}(z)&=\frac{ie^{-iz-b_{1}}+\frac{1}{i}e^{iz+b_{1}}}{2}\\
        &=\sin(z-b_{1}i).
\end{align*}
From \eqref{E:6.4} and \eqref{E:6.67b}, we get
\begin{align*}
w_{2}(-iz+c)&=\frac{e^{iz+b_{1}}-e^{-iz-b_{1}}}{2i}\\
            &=\sin(z-b_{1}i),\\
w_{2}(z)&=\frac{e^{i(iz-ic)+b_{1}}-e^{-i(iz-ic)-b_{1}}}{2i}\\
        &=\frac{e^{-z+b_{2}+(b_{1}-b_{2}+c)}-e^{z-b_{2}-(b_{1}-b_{2}+c)}}{2i}\\
        &=\frac{ie^{-z+b_{2}}-\frac{1}{i}e^{z-b_{2}}}{2i}\\
        &=\cosh{(z-b_{2})}.
\end{align*}
Clearly, $(w_{1},w_{2})=(\sin(z-b_{1}i),\cosh(z-b_{2}))$ is not a solution of Eq.\eqref{E:3.4}.

\par (2) If $a=-i$, we can set $p(z)=-iz+b_{1}$ and $h(z)=z+b_{2}$. From \eqref{E:6.66}, \eqref{E:6.67}, we get \begin{align}
\label{E:6.67c} e^{b_{1}+b_{2}-ic}&=1,\\
\label{E:6.67d} e^{b_{1}-b_{2}-c}&=-i.
\end{align}
Thus, by \eqref{E:6.68a} and \eqref{E:6.69}, we get
\begin{align*}
 w_{1}(-iz+c)&=\frac{e^{z+b_{2}}-e^{-z-b_{2}}}{2i}\\
             &=-i\sinh(z+b_{2}),\\
 w_{1}(z)&=\frac{-ie^{iz-b_{1}}+\frac{1}{-i}e^{-iz+b_{1}}}{2}\\
         &=\sin{(z+ib_{1})}.
\end{align*}
From \eqref{E:6.4} and \eqref{E:6.67d}, we gain
\begin{align*}
w_{2}(-iz+c)&=\frac{e^{-iz+b_{1}}-e^{iz-b_{1}}}{2i}\\
            &=-\sin(z+b_{1}i),\\
w_{2}(z)&=\frac{e^{-i(iz-ic)+b_{1}}-e^{i(iz-ic)-b_{1}}}{2i}\\
        &=\frac{e^{z+b_{2}+(b_{1}-b_{2}-c)}-e^{-z-b_{2}-(b_{1}-b_{2}-c)}}{2i}\\
        &=\frac{-i e^{z+b_{2}}-\frac{1}{-i}e^{-z-b_{2}}}{2i}\\
        &=-\cosh{(z+b_{2})}.
\end{align*}
We can also verify that $(w_{1},w_{2})=(\sin(z+b_{1}i),-\cosh(z+b_{2}))$ is not a solution of Eq.\eqref{E:3.4}.

\par Therefore,we have concluded that Eq.\eqref{E:3.4} does not posses any finite order of transcendental entire solutions other than $q=\pm 1$. Hence, our theorem follows. $\hfill\Box$









\end{document}